\documentclass[a4paper,12pt]{article}
\usepackage{amsmath}
\usepackage[english]{babel} 

\usepackage{graphicx,color} 
\DeclareGraphicsExtensions{.eps .pdf} 
\usepackage{epstopdf}
\usepackage{authblk}
\usepackage{float}
\usepackage{fixmath}
\usepackage{bm}
\usepackage{amsmath}

\usepackage{inputenc,t1enc}

\usepackage{latexsym}
\DeclareSymbolFont{AMSb}{U}{msb}{m}{n}
\DeclareSymbolFontAlphabet{\mathbb}{AMSb}

\begin{document}  

\bibliographystyle{plain}
 
\newcommand{\nc}{\newcommand}
\nc{\nt}{\newtheorem}
\nt{defn}{Definition}
\nt{lem}{Lemma}
\nt{pr}{Proposition}
\nt{theorem}{Theorem}
\nt{cor}{Corollary}
\nt{ex}{Example}
\nt{ass}{Assumption}
\nt{step}{Step}
\nt{case}{Case}
\nt{subcase}{Subcase}
\nt{note}{Note}
\nc{\bd}{\begin{defn}} \nc{\ed}{\end{defn}}
\nc{\blem}{\begin{lem}} \nc{\elem}{\end{lem}}
\nc{\bpr}{\begin{pr}} \nc{\epr}{\end{pr}}
\nc{\bth}{\begin{theorem}} \nc{\eth}{\end{theorem}}
\nc{\bcor}{\begin{cor}} \nc{\ecor}{\end{cor}}
\nc{\bex}{\begin{ex}}  \nc{\eex}{\end{ex}}
\nc{\bass}{\begin{ass}}  \nc{\eass}{\end{ass}}
\nc{\bstep}{\begin{step}}  \nc{\estep}{\end{step}}
\nc{\bcase}{\begin{case}}  \nc{\ecase}{\end{case}}
\nc{\bsubcase}{\begin{subcase}}  \nc{\esubcase}{\end{subcase}}
\nc{\bnote}{\begin{note}}  \nc{\enote}{\end{note}}
\nc{\prf}{{\it Proof} }
\nc{\eop}{\hfill $\Box$ \\ \\}
\nc{\argmin}{\mathrm{argmin}}
\nc{\argmax}{\mathrm{argmax}}
\nc{\sgn}{\mathrm{sgn}}
\nc{\Var}{\mathrm{Var}}
\nc{\Cov}{\mathrm{Cov}}
\nc{\bak}{\!\!\!\!\!}
\nc{\IBD}{\mathrm{IBD}}
\nc{\supp}{\mathrm{supp}}
\nc{\dom}{\mathrm{dom}}
\nc{\R}{{\mathbb R}}
\nc{\peq}{\preceq}
\nc{\wt}{\widetilde}
\nc{\Mult}{\mathrm{Mult}}
\newcommand{\tcb}{\textcolor{black}}
\newcommand{\tcr}{\textcolor{black}}
\nc{\Prob}[1]{\mathbb{P}_{#1}}
\newcommand{\notex}[1]
{$^{(!)}$\marginpar[{\hfill\tiny{\sf{#1}}}]{\tiny{\sf{(!) #1}}}}
\newcommand\approxin
  {\raisebox{-1ex}{ $\stackrel{\textstyle\in}{\scriptstyle\sim}$ }}

\title{A stochastic process approach to multilayer neutron detectors}
\author[1]{Dragi Anevski}
\author[2]{Richard Hall-Wilton}
\author[2]{Kalliopi Kanaki}
\author[1]{Vladimir Pastukhov \thanks{pastuhov@maths.lth.se} }


\affil[1]{Centre for Mathematical Sciences, Lund University, Lund SE-221 00, Sweden}
\affil[2]{European Spallation Source ERIC, P.O Box 176, Lund  \newline SE-221 00, Sweden}
\date{}
\maketitle
\begin{abstract}

The sparsity of the isotope helium-3, ongoing since 2009, has initiated a new generation of neutron detectors. One particularly promising development line for detectors is the multilayer gaseous detector. 
In this paper, a stochastic process approach is used to determine the neutron energy from the additional data afforded by the multilayer nature of these novel detectors. 

The data from a multilayer detector consists of counts of the number of absorbed neutrons along the sequence of the detector's layers, in which the neutron absorption probability is unknown. We study the maximum likelihood estimator for the intensity and absorption probability, and show its consistency and asymptotic normality, as the number of incoming neutrons goes to infinity. We combine these results with known results on the relation between the absorption probability and the wavelength to derive an estimator of the wavelength and to show its consistency and asymptotic normality. 

\end{abstract}

\section{Introduction}
The European Spallation Source\footnote{\tt https://europeanspallationsource.se} (ESS), sited in Lund,
Sweden, is planned to be operational in
2019 and the world's leading source for the study of materials using
neutrons by 2025.

In order to address the challenge  of developing a new generation
of neutron detectors an international collaboration of 10 neutron scattering institutes in
Europe, Asia and America  (the International Collaboration on the Development of
  Neutron Detectors\footnote{{\tt http://icnd.org}}) was formed in 2010. The members have chosen as the three most promising technologies for investigation: Scintillator detectors,
boron-10 thin film detectors and $^{10}$BF$_{3}$ gas detectors. At present boron-10 thin film detectors seem to be the only realistic solution
for large area detectors ($>$ 10 m$^2$ active detector area).  For the ESS, novel neutron detectors represent a critical technology that needs to be developed, with corresponding research and development done as contributions to the ESS design work.  

In this paper we study the feasibility and possibility of the statistical determination of neutron wavelength for the new generation of neutron detectors being developed at the ESS.

Assume that a beam of neutrons arrives at the face of the detector. The detector consists of a sequence of boron-10 coated layers, between which there are gas-filled cavities. The principle of the detector can be described in a simplified manner as follows. A neutron that goes through a boron-10 layer can sometimes interact with a boron-10 atom in the layer, temporarily exciting the atom into an unstable state from which it will fall back to a stable state and thereby emit an electrically charged particle, that will ionise the gas. This electrical potential in the gas filled chamber is detected and the instrument notes that a neutron has been absorbed, see  Kanaki  \textit{et al.} (2013).  The outcome of this is that we have a count of +1 in the number of neutrons that have passed and been detected.
The probability with which a neutron is absorbed and detected is a function of the energy content of the neutron, i.e.\ a function of the neutron wavelength.

If we view the neutron beam as a set of particles that hit the face of the detector, then each neutron will either be absorbed or not at the first layer. If the neutron is not absorbed at the first layer, it may possibly be absorbed at the second layer, and so on. From the simplified description above it is clear that data from a multilayer detector will consist of counts of the number of absorbed neutrons along the sequence of the detector's layers.

By a beam we mean a stream of particles with a certain fixed wavelength $\mu$. Let the number of neutrons that arrive in the time interval $[0,t]$ be denoted by $X_0(t)$. Then $X_0(t)$ is a counting process, such that $X_0(0)=0$. 

A simple model for the process of incoming neutrons $X_0(t)$ is that of a Poisson process with intensity $\lambda$. The Poisson model assumption is reasonable since neutrons are electrically neutral particles and there are therefore no long-distance interactions between the particles in the beam, see Willie \& Carlile (1999), Chapter 2, for a discussion of the model. 
The intensity $\lambda$ is assumed to be an unknown nuisance parameter, and will be estimated.

At a layer each neutron is absorbed with a certain probability $p$ (the absorption efficiency). The probability of absorption $p$ is also assumed to be an unknown parameter,  its dependence on the wavelength $\mu$ of the incident neutron is, however, of a known functional form, see Kanaki  \textit{et al.} (2013). This property will be used to make  inference about the parameter $\mu$. For a more thorough introduction to the subject of neutron interactions we refer to Chapter 2 in Willie \& Carlile (1999).

As will be shown later, our data set is generated by a sequentially thinned Poisson process, which is a special case of multinomial thinning. Inference for thinned point processes was studied in detail in Karr (1985) and Bensa{\"i}d (1997), where, in particular, the authors studied the problem of estimation of the thinning parameter $p$ from observation of the thinned processes. In Karr (1985) and Bensa{\"i}d (1997) the thinning parameter $p$ is defined as a function from an underlying compact metric space to $[0,1]$. In Karr (1985) the author uses a nonparametric histogram estimator of $p$ and in Bensa{\"i}d (1997) the author studies a kernel estimator. 

Though the approaches developed in Karr (1985) and Bensa{\"i}d (1997) are quite general, they cannot be applied to the problem considered in this paper because, first, in our case the absorption probability (thinning parameter) is homogeneous (does not depend on the time of experiment) and, therefore, we can use a parametric approach to estimate it and, second, our data come from a multinomial thinning of the original Poisson process, not a binomial one as in Karr (1985) and Bensa{\"i}d (1997). 

The problem of multinomial thinning of point processes was studied in Long (1995), where, in particular, the author proved that a point process is Poisson if and only if the thinned processes are independent and Poisson. However,  to our knowledge, the problem of inference for a sequentially thinned Poisson process has not been studied yet.
Given the data, we suggest in this paper a likelihood approach and study the maximum likelihood estimator (mle) of the two-dimensional parameter $(\lambda,p)$, where $\lambda$ is the intensity and $p$ the thinning parameter (absorption probability). We derive conditions for the existence of the mle and prove its consistency and asymptotic normality, as the experiment time (or number of incoming neutrons) goes to infinity. We combine these results with known results for the relation between the absorption probability and the wavelength to derive a final estimator of the wavelength and to show consistency and asymptotic normality for the estimator. We also state results on the precision of the estimator, by deriving a relation between the width of the confidence interval, for the unknown wavelength, and the detector construction, in  terms of  the number of layers used in the detector. The performance of the estimator is illustrated on simulated data.

There are two main results of this paper. 
The first establishes the feasibility of estimating the wavelength of a neutron beam, based only on count data of the number of detected neutrons. The second determines necessary features of the detector, which for the specific detector is the number of layers, in order to be able to estimate the wavelength with a given precision. 
Following the construction of the ESS research facility, we intend to apply our estimation procedures to experimental data.

The paper is organized as follows. Section \ref{secstocmod} provides the general scheme of the neutron detector and the modeling of neutron interactions with the detector layers. Section \ref{secinfpar} is devoted to inference of the parameters. We derive the mle for the intensity $\lambda$ of an incident beam  and absorption efficiency $p$, in Lemma \ref{Lmzero} and \ref{Aoneroot} we discuss the uniqueness of the solutions to the score equations, and in Theorem \ref{consnormef}, which is one of the main results of this paper, we derive the strong consistency and asymptotic normality of the mle. In Corollaries \ref{cormunorm} and \ref{cormuest} we derive the consistency and asymptotic normality of the mle of the wavelength. Using these final results we are able to construct confidence intervals for the wavelength. Section \ref{secsimexp} gives a simulation study to explore the estimator's performance. Section \ref{secdisc} contains a discussion of the results presented in the paper and plans for future work. Proofs of all results are given in the Appendix.

\section{Scheme of a discrete spacing detector}\label{secstocmod}

\begin{figure}[h]
\center
\includegraphics[width=1\linewidth]{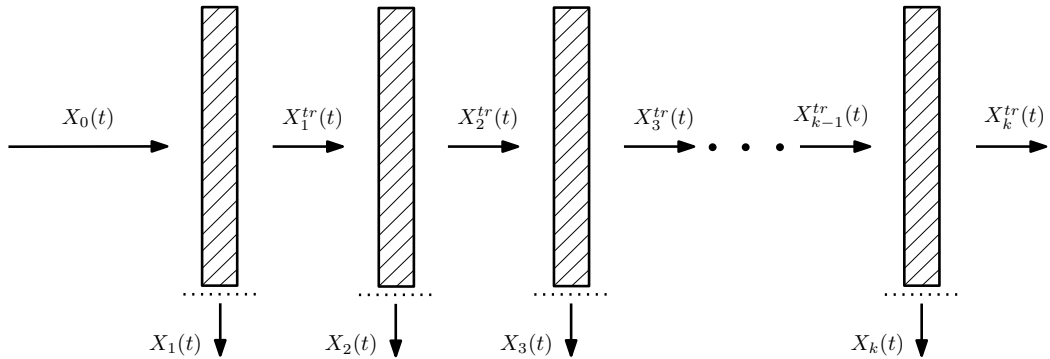}
\caption{\label{detektork} The scheme of the detector.}
\end{figure}

Assume that an incident beam of neutrons hits the first layer of the detector, cf. Fig.\ref{detektork}. At the layer a neutron can possibly be absorbed and detected. If a neutron is not absorbed it will go through the detector's layer. We assume that these are the only two possibilities for the neutron interaction with a layer, i.e.\ it is assumed that the probability of an inelastic scattering of a neutron in the boron layers or in the material of the layers is negligibly small. Let $p$ be the probability of absorption of a neutron, so that $1-p$ is the probability of its transmission. If a neutron is absorbed, it will then be detected. Let $X_{1}(t)$ be the number of neutrons that are absorbed at the first layer, so that $X^{tr}_{1}(t) =X_0(t)-X_{1}(t)$ is the number of transmitted neutrons.

Now assume that the beam of transmitted neutrons $X^{tr}_{1}(t)$ hits the second layer, at which, again, each neutron can either be absorbed (with the same probability $p$ as at the previous layer) and  then detected, or transmitted to the second layer. Let $X_{2}(t)$ be the number of neutrons that are absorbed at the second layer and let $X^{tr}_{2}(t)=X^{tr}_{1}(t) - X_{1}(t)$  be the number of transmitted neutrons. We assume that the registrations (absorptions) of different particles are independent and the times of absorption and travelling from layer to layer are negligibly small. This behaviour is repeated at each layer and gives the general scheme for the neutron beam's absorption and transmission in the detector. 

Let $X_{i}(t)$ be the number of neutrons absorbed at the layer $i$ in the time interval $[0, t]$ and let $X^{tr}_{i}(t)$ be the number of transmitted neutrons in the same time interval through the layer $i$, for $i=1, \dots, k$. Then $X_{i}(t)$ and $X^{tr}_{i}(t)$ are counting processes and $X_{i}(0) = 0$ and $X^{tr}_{i}(0) = 0$, for $i=1, \dots, k$. 
\blem\label{indepcounts}

The processes $\{X_{i}(t)\}_{i\geq 1}$ are jointly independent Poisson processes with intensities $p(1-p)^{i-1}\lambda$, respectively.
\elem
The statement of Lemma \ref{indepcounts} follows from the property of a multinomial thinning of a Poisson process cf. Theorem 5.17 in   Kulkarni (2009), Long (1995), Assuncao R. M. \& Ferrari P. A. (2007).

\section{Inference for the parameters}\label{secinfpar}
Now suppose that we have run an experiment at the neutron detector, the result of which is a sequence of counts of the numbers of detected neutrons along the detector. Let us denote the data as a vector $x=(x_1,\ldots,x_k)$ of integers, with $x_i$ being the number of observed neutrons at layer $i$, for $i=1,\ldots,k$. From Lemma \ref{indepcounts} we know that the data are observations of independent Poisson distributed random variables, with unknown expectations $p(1-p)^{i-1}\lambda$, for $i=1, \dots, k$. 


\subsection{The mle of the thinning parameter $p$ and the intensity of an incident process $\lambda$}
We are interested in deriving consistency and asymptotic normality of the estimators. For this we need to explain what we mean by letting "the amount of data" go to infinity. There are several ways to model this. We can either let the experiment time $t$ increase, or we can view the problem as a repeated measurement problem and thus make several, $n$ of them, independent measurements during a fixed time interval $[0,t]$ and instead let $n$ go to infinity. 
Since we use the Poisson process as a model for the neutron beam, the two approaches will give quantitatively the same limit results. We choose to view the problem as a repeated sample problem. 

The inference problem can be described as follows. We perform $n$ experiments. For each experiment $j=1,\ldots,n$, we measure the number of neutrons $X_{ij}$ detected at layer $i=1\ldots,k$ during the time interval $[0,t]$. Thus  $\{X_{ij}\}_{i,j=1}^n$ are random variables and $\{x_{ij}\}_{i,j=1}^n$ are the values which they take. Let $(p, \lambda)$ denote the parameters, that are assumed to lie in $[0,1]\times [0,\infty)$. Introduce the vectors
$\mathbf{X}_{j}=(X_{1j},\ldots,X_{kj})^{T}$ and $\mathbf{x}_{j}=(x_{1j},\ldots,x_{kj})^{T}$, respectively. Note that  the vectors $\mathbf{X}_{j}$ are independent random vectors with mutually independent components $X_{ij}$, by Lemma \ref{indepcounts}, from $n$ independent experiment rounds. Finally denote $\mathbf{X} = [\mathbf{X}_{1},\ldots, \mathbf{X}_{n}]$ and $\mathbf{x} = [\mathbf{x}_{1},\ldots, \mathbf{x}_{n}]$, and note that these are $k\times n$ matrices of discrete random variables and of integers values, respectively. 

Thus we let $X_{ij}$ be the number of neutrons observed at the layer $i$ at the experiment round $j$ with probability mass function
\begin{eqnarray*}
f(x_{ij} \vert p, \lambda) = e^{-m_i}\frac{m_i^{x_{ij}}}{x_{ij}!},
\end{eqnarray*}
where $m_{i} = p(1-p)^{i-1}\lambda t$. Then each vector $\mathbf{X_j} = (X_{1j},\ldots,X_{kj})^{T}$ has the joint distribution 
\begin{eqnarray*}
f(\mathbf{x}_{j} \vert p, \lambda)= \prod_{i=1}^k f(x_{ij} \vert p, \lambda) = \prod_{i=1}^k e^{-m_i}\frac{m_i^{x_{ij}}}{x_{ij}!}.
\end{eqnarray*}

Note, that if $k=1$, then $m=p\lambda t$ and, therefore, in this case one can only estimate the product $p\lambda$, and not $p$ and $\lambda$  separately.

Assume that $k > 1$. The log-likelihood is then given by
\begin{eqnarray}\label{loglikehn}
     l_{n}(p, \lambda \vert \mathbf{x})
     &=&\sum_{j=1}^n\sum_{i=1}^k( -  m_i + x_{ij}\log m_i - \log x_{ij}! )\nonumber.
 \end{eqnarray}
The mle $(\hat{p}_{n},\hat{\lambda}_{n})$  is the solution of the score equations
\begin{eqnarray}\label{psilikk}
&&\left\{ 
\begin{array}{ll}
&\frac{1}{n}\frac{\partial l_n}{\partial \lambda} = \frac{s_n-\lambda t(1-(1-p)^k)}{\lambda} = 0,\\
&\frac{1}{n}\frac{\partial l_n}{\partial p} = \frac{(1-p)(s_n+z_n)-z_n-\lambda t(k(1-p)^k-k(1-p)^{k+1})}{p(1-p)} =0,
\end{array}
\right.
\end{eqnarray}
where $s_n=\frac{1}{n}\sum_{j=1}^n \sum_{i=1}^k x_{ij}$ and $z_n=\frac{1}{n}\sum_{j=1}^n\sum_{i=1}^k (i-1)x_{ij}$.
If we assume that $\hat{p}_{n}(1-\hat{p}_{n})\neq 0$, $\hat{\lambda}_{n}\neq 0$ we get the system of equations
\begin{eqnarray}
&&\left\{ 
\begin{array}{ll}
&s_n-\hat{\lambda}_{n}t(1-\hat{y}_{n}^k) = 0,\\
&a_{n}\hat{y}_{n}^{k+1} - b_{n} \hat{y}_{n}^k + c_{n} \hat{y}_{n} - d_{n} = 0, \label{eq:yn-eq}
\end{array}
\right.
\end{eqnarray}
where 
\begin{eqnarray}\label{abcdx}
a_n&=& - s_n - z_n + ks_n,\nonumber\\
b_n&=& - z_n + ks_n,\nonumber\\
c_n&=&z_n + s_n,\\
d_n&=&z_n,\nonumber\\
\hat{y}_{n}&=&1-\hat{p}_{n}.\nonumber
\end{eqnarray}
Obviously $(\ref{eq:yn-eq})$ has exactly one solution $(\hat{p}_{n},\hat{\lambda}_{n})$ if and only if the second equation in $(\ref{eq:yn-eq})$ has exactly one root. 

\blem \label{Lmzero}
The function 
\begin{eqnarray*}
f(y)=a_{n}y^{k+1} - b_{n} y^k + c_{n} y - d_{n},
\end{eqnarray*}
for $k>1$ with coefficients given in (\ref{abcdx}),
has one zero in the open interval $(0,1)$ when the inflection point $y_{i.p.}$ satisfies the inequality 
\begin{eqnarray*}
    y_{i.p.} := \frac{b_n(k-1)}{a_n(k+1)}<1, \label{eq:x0def}
\end{eqnarray*}
and no zeros in $(0,1)$ when $y_{i.p.}\geq 1$. 
\elem

Lemma \ref{Lmzero} gives the condition of existence and uniqueness of $(\hat{p}_{n},\hat{\lambda}_{n})$,
but there is no guarantee that it holds for a finite $n$. However, the following result holds.

\blem\label{Aoneroot}
Let $A_n=\{\mbox{Equation (\ref{eq:yn-eq}) has exactly one root in $(0,1)$} \}$. Then $A_{n}$ happens for all sufficiently large $n$ almost surely.
\elem

\subsubsection{Asymptotic properties of the mle}

\bth\label{consnormef}

The mle $(\hat{p}_n, \hat{\lambda}_n)$,  given in (\ref{psilikk}), is strongly consistent
\begin{eqnarray*}
(\hat{p}_n, \hat{\lambda}_n) \stackrel{a.s.}{\to} (p, \lambda),
\end{eqnarray*}
and asymptotically normal
\[ 
 \sqrt{n}((\hat{p}_n, \hat{\lambda}_n) - (p, \lambda)) \stackrel{d}{\to}  \mathcal{N}( \mathbf{0}, \left[\mathbf{I}(p, \lambda)\right]^{-1} ),          
\]
as $n \to \infty$, where $\mathbf{I}(p, \lambda)$ is the information matrix

\begin{eqnarray*}
\mathbf{I}(p, \lambda) &=& \frac{1}{k}\sum_{i=1}^k \mathbf{I}_{(i)}(p, \lambda), \label{eq:info.matrix}
\end{eqnarray*}
where $\mathbf{I}_{(i)}(p, \lambda)$ denotes the information matrix corresponding to $f(x_{ij} \vert p, \lambda)$ with fixed $i$.
\eth

From the theorem above, after simplification, we obtain the following asymptotic covariances
\begin{eqnarray*}\label{elmIsigma}
    \sigma^{2}_{p}(p, \lambda) = \left[\mathbf{I}(p, \lambda)\right]^{-1}_{p p} &=& \frac{(1-(1-p)^{k})(1-p)p^{2}}{\lambda t q(p,k)}\to\frac{(1-p)p^{2}}{\lambda t},\nonumber \\
     \sigma^{2}_{\lambda}(p, \lambda) = \left[\mathbf{I}(p, \lambda)\right]^{-1}_{\lambda \lambda} &=& \frac{\lambda h(p,k)}{t q(p,k)}\to \frac{\lambda}{t}, \\
    \sigma^{2}_{p, \lambda}(p, \lambda) = \left[\mathbf{I}(p, \lambda)\right]^{-1}_{\lambda p} &=& \frac{kp((1 - p)^{k} - (1 - p)^{k-1})}{t q(p,k)}\to 0,\nonumber
\end{eqnarray*} 
as $k\to\infty$,
where
\begin{eqnarray*}
h(p,k) = 1-k^{2}(1-p)^{k+1} + (2k^{2} - 1)(1-p)^{k} - k^{2}(1-p)^{k-1},
\end{eqnarray*}
and
\begin{eqnarray}\label{eq:q(p,k)-def}
q(p,k) &=& (1-p)^{2k} - k^{2}(1-p)^{k+1} + 2(k^{2} - 1)(1-p)^{k}\nonumber\\
&& - k^{2}(1-p)^{k-1} + 1. 
\end{eqnarray}

We are mainly interested in the estimation of $p$, since there is a functional relation between the absorption efficiency $p$ and  the wavelength $\mu$ of the incident neutrons, cf. $(\ref{prob1})$ and $(\ref{eq:prob2})$ below.  Analysing the behaviour of $\sigma^{2}_{p}(p, \lambda)$, it can be shown that $\sigma^{2}_{p}(p, \lambda)$ is a strictly decreasing function of $k$ for every $p \in (0, 1)$.





\subsection{Estimation of the wavelength $\mu$ of an incident beam.}


We are interested in estimating the wavelength of a monochromatic neutron beam. The probability of absorption $p$ depends on the neutron wavelength $\mu$ as (cf. Willie \& Carlile (1999), Section 2.3)
\begin{eqnarray}\label{prob1}
p = 1 - e^{-\Sigma(\mu)\rho_{at}d_{l}},
\end{eqnarray} 
where the parameter $\Sigma(\mu)$ is called the cross-section of absorption, $\rho_{at}$ is the atomic density of $^{10}B$ in the $B_{4}C$ coating and  $d_{l}$ is the thickness of the boron layer. Example values of parameters in a detector are $\rho_{at} = 10^{29}$ $m^{-3}$, $d_{l} = 10^{-6}$ $m$, cf.  Kanaki \textit{et al.} (2013). 

The neutron cross-section $\Sigma(\mu)$ can be modelled as 
\begin{eqnarray*}
\Sigma(\mu) = \varsigma\mu, \label{eq:prob2}
\end{eqnarray*} 
where the coefficient $\varsigma$ is different for different materials, see Willie \& Carlile (1999), cf. Section 2.3. Furthermore, the coefficient  $\varsigma$ does not depend on the neutron wavelength and has been measured experimentally, cf. Schmitt \textit{et al.} (1959). From the results in Schmitt \textit{et al.} (1959) we conclude that the estimator $\hat{\varsigma}$ of $\varsigma$ is unbiased and asymptotically normal
\begin{eqnarray*}
\sqrt{n'}(\hat{\varsigma}_{n'} - \varsigma)   \stackrel{d}{\to} \mathcal{N}(0,  \sigma^{2}_{\varsigma}),
\end{eqnarray*}
as $n' \to\infty$, where $n'$ is the number of runs performed in the experiment to estimate $\varsigma$ and $\sigma^{2}_{\varsigma}$ is its asymptotic variance. 

Let us rewrite (\ref{prob1}) as
\begin{eqnarray}\label{pmu}
p = 1 - e^{-\chi\mu},
\end{eqnarray} 
where 

\begin{eqnarray}\label{eq:chi-varsigma}
\chi = \rho_{at}d_{l}\varsigma. 
\end{eqnarray}
Then, from delta method,  the plug-in estimator $\hat{\chi}=\rho_{at}d_{l}\hat{\varsigma}$ of $\chi$ is asymptotically normal 
\begin{eqnarray}
\sqrt{n'}(\hat{\chi}_{n'} - \chi) \stackrel{d}{\to} \mathcal{N}(0, \sigma^{2}_{\chi}), \label{eq:xin_n-est}
\end{eqnarray} 
with $\chi = \rho_{at}d_{l}\varsigma$ and $\sigma^{2}_{\chi} = \rho^{2}_{at}d^{2}_{l}\sigma^{2}_{\varsigma}$.

From (\ref{pmu}), we obtain
\begin{eqnarray}\label{muonphi}
\mu (p, \chi) = -\frac{\log(1-p)}{\chi}.
\end{eqnarray} 

Next, we combine two limit distribution results, for $\hat{p}_n$ and for $\hat{\chi}_{n'}$, to get a limit distribution for the plug-in estimator of $\mu$. In order to formalize this in a proper way, we introduce a factor $\gamma$, which is merely the (asymptotic) ratio between $n'$ and $n$. The result in a practical finite-sample situation will be used in exactly that way: by letting $\gamma=n'/n$ and using the limit distribution to provide asymptotic confidence intervals or tests.

\bcor\label{cormunorm}
The plug-in estimator $\hat{\mu} = \mu (\hat{p}_{n}, \hat{\chi}_{n'})$ of $\mu$ is asymptotically normal
\begin{eqnarray*}\label{cltmu}
\sqrt{n}(\hat{\mu} - \mu) \stackrel{d}{\to} \mathcal{N}(0, \sigma^{2}_{\mu}), 
\end{eqnarray*}
as $n\to\infty$, with 
\begin{eqnarray*}
\sigma^{2}_{\mu} = \left[ \frac{\partial\mu}{\partial p}(p, \chi)\right]^{2}\sigma^{2}_{p}(p, \lambda) +  \frac{1}{\gamma} \left[\frac{\partial\mu}{\partial\chi}(p, \chi)\right]^{2}\sigma^{2}_{\chi},
\end{eqnarray*}
where $n$ is the number of measurements for $\hat{p}_{n}$ and $n' = \lceil \gamma n \rceil$, $\gamma  > 0$, is the number of measurements for  $\hat{\chi}_{n'}$ ($\lceil \gamma n \rceil$ is the smallest integer not less than $\gamma n$).
\ecor


Introduce the notation
\begin{eqnarray}\label{estsmu}
S^{2}_{n}(\hat{p}, \hat{\lambda}, \hat{\chi}) = \left[ \frac{\partial\mu}{\partial p}(\hat{p}_{n}, \hat{\chi}_{n'})\right]^{2}\sigma^{2}_{p}(\hat{p}_{n}, \hat{\lambda}_{n}) + \frac{1}{\gamma} \left[\frac{\partial\mu}{\partial\chi}(\hat{p}_{n}, \hat{\chi}_{n'})\right]^{2}\hat{\sigma}^{2}_{\chi}, 
\end{eqnarray}
where both the estimate $\hat{\chi}_{n'}$ and the estimate of the variance $\hat{\sigma}^{2}_{\chi}$ are based on $n'$ measurements, and $(\hat{p}_{n}, \hat{\lambda}_{n})$ are the mle of $(p, \lambda)$ based on $n$ measurements. 

The next result follows from Slutsky's theorem and the continuous mapping theorem, cf. Chapter 2 in  van der Vaart (1998).
\bcor\label{cormuest}
Under the assumptions of the previous Corollary
\begin{eqnarray*}\label{confint}
\frac{\sqrt{n}(\hat{\mu} - \mu )}{S_{n}} \stackrel{d}{\to} \mathcal{N}(0, 1), 
\end{eqnarray*}
as $n \to \infty$.
\ecor

Using the above limit distribution result for the mle $\hat{\mu}$ we can construct an approximate $100(1 - \alpha)$ per cent  confidence interval for $\mu$, viz.:
\begin{eqnarray}\label{confmu}
[\mu (\hat{p}_{n}, \hat{\chi}_{n'}) - z_{\alpha/2} \frac{S_{n}}{\sqrt{n}},\ \mu (\hat{p}_{n}, \hat{\chi}_{n'}) + z_{\alpha/2} \frac{S_{n}}{\sqrt{n}} ],
\end{eqnarray}
where $z_{\alpha/2}$ is the $\alpha/2$-th quantile of the standard normal distribution.

Let us rewrite the expression for $\frac{(S_{n})^{2}}{n}$ as

\begin{eqnarray}\label{Sn}
\frac{(S_{n})^{2}}{n} =  (S^{(p)}_{\hat{\mu}})^{2} + (S^{(\chi)}_{\hat{\mu}})^{2}, 
\end{eqnarray}
where 
\begin{eqnarray}\label{Sp}
S^{(p)}_{\hat{\mu}}(\hat{p}_{n}, \hat{\lambda}_{n}, \hat{\chi}_{n'}) &=& \frac{1}{\sqrt{n}} \frac{\partial\mu}{\partial p}(\hat{p}_{n}, \hat{\chi}_{n'})\sigma_{p}(\hat{p}_{n}, \hat{\lambda}_{n})
=\frac{\sigma_{p}(\hat{p}_{n}, \hat{\lambda}_{n})}{\sqrt{n}(1 - \hat{p}_{n})\hat{\chi}_{n'}},
\end{eqnarray} 
\begin{eqnarray}\label{Shi}
S^{(\chi)}_{\hat{\mu}}(\hat{p}_{n}, \hat{\chi}_{n'}) &=&\frac{1}{\sqrt{\gamma n}}\frac{\partial\mu}{\partial\chi}(\hat{p}_{n}, \hat{\chi}_{n'}) \hat{\sigma}_{\chi} = \frac{1}{\sqrt{\gamma n}}\frac{\log(1 - \hat{p}_{n})}{\hat{\chi}_{n'}^{2}}\hat{\sigma}_{\chi}.
\end{eqnarray} 
Next, since $\gamma$ is the asymptotic ratio between $n'$ and $n$, then one can rewrite (\ref{Shi}) as
\begin{eqnarray}\label{Shi1}
S^{(\chi)}_{\hat{\mu}}(\hat{p}_{n}, \hat{\chi}_{n'}) &\approx& \frac{\log(1 - \hat{p}_{n})}{\sqrt{n'} \hat{\chi}_{n'}^{2}}\hat{\sigma}_{\chi}
\end{eqnarray}
for relatively big values of $n'$. 

The coefficient $\gamma$ takes into account that the number of experimental runs $n$ for estimating $p$ and $\lambda$ is not equal to $n'$, which is the size of the sample used in the estimation of $\chi$.  We emphasise, that in the simulation experiments belows the value of $n'$ is fixed and, therefore, $S^{(\chi)}_{\hat{\mu}}$ does not decrease with increasing $n$. Therefore, we can view this term as a kind of systematic error, outside of our control.

\section{A simulation experiment}\label{secsimexp}

\begin{figure}[t]
\center
\includegraphics[width=1\linewidth]{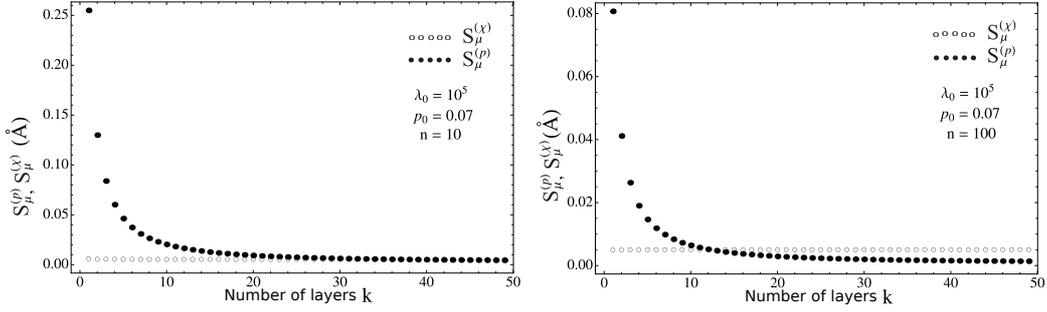}
\caption{\label{fci} The dependence of $S^{(p)}_{\hat{\mu}}$ and $S^{(\chi)}_{\hat{\mu}}$ on the number of layers $k$.}
\end{figure}

\begin{figure}[b]
\center
\includegraphics[width=1\linewidth]{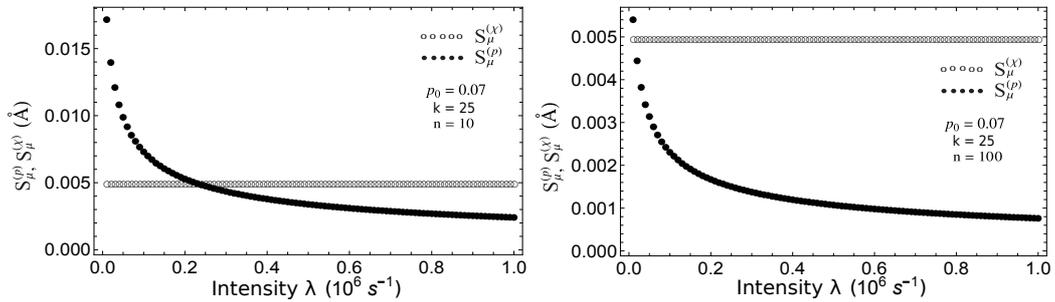}
\caption{\label{Sl} The dependence of $S^{(p)}_{\hat{\mu}}$ and $S^{(\chi)}_{\hat{\mu}}$ on the the intensity of an incident beam $\lambda$.}
\end{figure}

In this section we perform a simulation experiment to evaluate the estimator's performance. In particular, we illustrate the dependence of individual terms in (\ref{Sn}) on the number of layers (Figure \ref{fci}) and on the intensity of a beam (Figure \ref{Sl}), and the confidence interval width's dependence on the number of layers for several wavelengths (Figure \ref{muonn10}). 

We simulate a Poisson process $X_{0}(t)$ a number of times $n$,  for $n = 10, \ 100$, for the parameters values $p = 0.05, \ 0.07, \ 0.1$, $\lambda = 10^{5}$ $s^{-1}$, which correspond to the wavelengths $\mu$ = 2.4, 3.4 and 4.9 \AA. These are typical neutron wavelengths for the possible applications of the detector, see Kanaki \textit{et al.} (2013). 

The mle $(\hat{p}_{n},\hat{\lambda}_{n})$ is calculated for the simulated data.  We recall the relation between $\chi$ and $\varsigma$ in $(\ref{eq:chi-varsigma})$, and note that  $\rho_{at}$ and $d_{l}$ are known. The estimator of $\varsigma$ is assumed to be asymptotically normal, with mean value the sample mean and variance equal to a pooled variance estimate using three series of 15 measurements, which give in total $n'=45$ experimental data points, see Schmitt \textit{et al.} (1959). Using the results of Schmitt \textit{et al.} (1959) we have the following estimates for $\chi$: $\hat{\chi}_{n'} = 2.142 \times 10^8$ $m^{-1}$ and $\hat{\sigma}^{2}_{\chi} = 0.021 \times 10^8$ $m^{-2}$.

First, we analyse the dependence of the approximate confidence interval on the number of detector layers. 
Figure \ref{fci} shows the dependence of $S^{(p)}_{\hat{\mu}}$ and $S^{(\chi)}_{\hat{\mu}}$, defined in (\ref{Sp}) and (\ref{Shi}), on the number of the layers in the detector for 10 and 100 runs of the experiment. We note, in particular, that $S^{(p)}_{\hat{\mu}}$ and $S^{(\chi)}_{\hat{\mu}}$ are of the same size at $k \approx 25$ for $n = 10$ experimental runs and at $k \approx 15$ for $n = 100$.

Second, we study the dependence of the approximate confidence interval on the intensity $\lambda$ of an incident beam. Figure \ref{Sl} displays the dependence of $S^{(p)}_{\hat{\mu}}$ and $S^{(\chi)}_{\hat{\mu}}$ on $\lambda$ for 10 and 100 runs of the experiment for the fixed number of layers $k=25$. One can see that if $n=10$ the term $S^{(\chi)}_{\hat{\mu}}$ becomes dominant when $\lambda  > 0.2 \times 10^{6}$ and if $n=100$ the term $S^{(\chi)}_{\hat{\mu}}$ dominates $S^{(p)}_{\hat{\mu}}$ even for small intensities ($\lambda < 10^{5}$)

\begin{figure}[t]
\center
\includegraphics[width=0.7\linewidth]{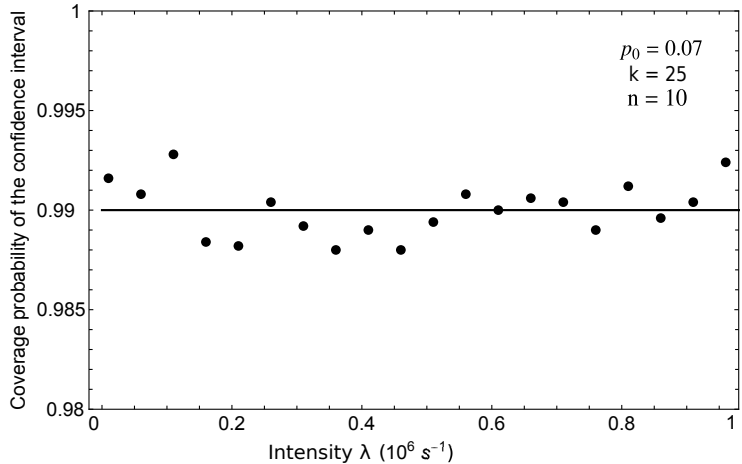}
\caption{\label{covprob} Dependence of the coverage probability of the approximate confidence interval on the intensity of an incident beam $\lambda$.}
\end{figure}

\begin{figure}[h]
\center
\includegraphics[width=0.9\linewidth]{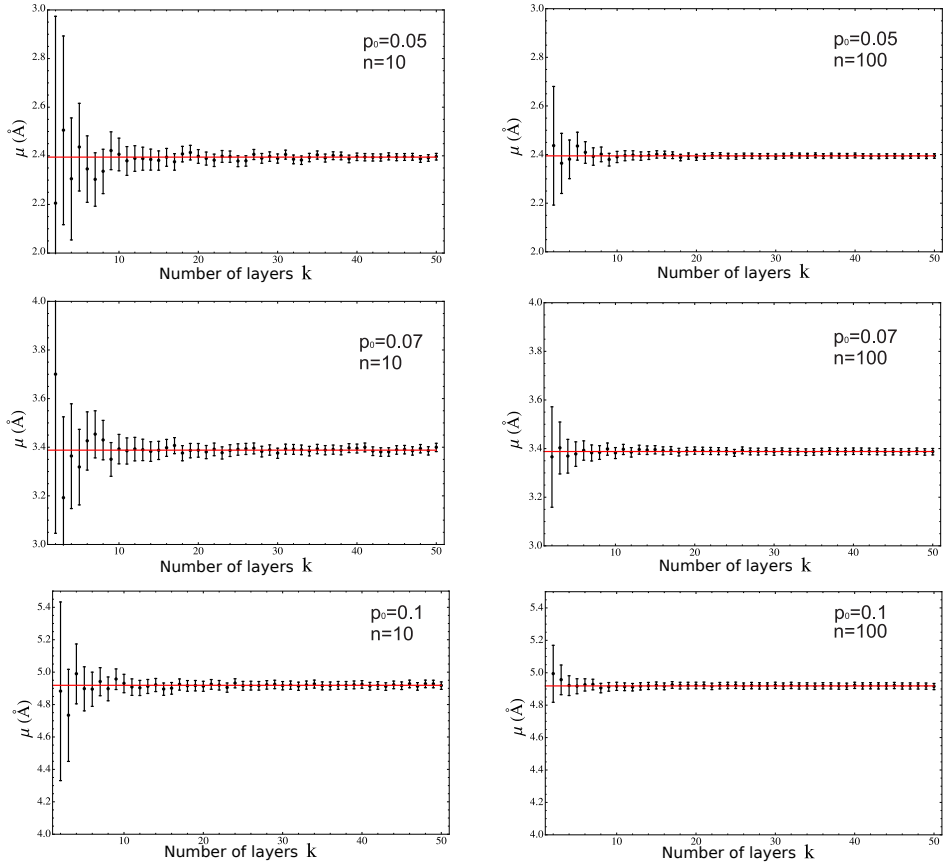}
\caption{\label{muonn10}99\% confidence interval for $\mu$ based on simulations for $n=10$, $100$ and $p=0.05$, $0.07$, $0.1$, $\lambda=10^{5}$ $s^{-1}$, $t=1$s. The red line is the true value of $\mu$.}
\end{figure}

Next, in order to assess the accuracy of the asymptotic approximation we estimate the coverage probability of the approximate confidence interval based on 5000 Monte-Carlo simulations. From Figure \ref{covprob} one can see that the deviation of the confidence band's coverage probability is less that 0.5 \% even for the quite small number of repetitions $n=10$. 

In Figure \ref{muonn10} we have plotted the confidence interval bars as a function of the number of layers, for $\mu$ = 2.4, 3.4 and 4.9 \AA \, and $n = 10$ and $100$.

The results of the simulation experiments show that the errors are rapidly decreasing as a function of the number of layers $k$ in the detector, cf. Figure 2, where the term $S^{(p)}_{\hat{\mu}}$ we may control by increasing the number of measurements, whereas the term $S^{(\chi)}_{\hat{\mu}}$ we are not able to influence and therefore we can see as a form of systematic error contribution to the total variance (\ref{Sn}). As indicated in Figure 2, for the choice of model parameters, at approximately 10-25 layers the term $S^{(p)}_{\hat{\mu}}$ that we can affect becomes smaller than the systematic error term $S^{(\chi)}_{\hat{\mu}}$. Figure 3 shows that, again, the term $S^{(p)}_{\hat{\mu}}$ decreases with increasing intensity, whereas the term $S^{(\chi)}_{\hat{\mu}}$ is almost not affected by a change in intensity.

Note that in our simulations for Figure 4, and only here, in our assessment of the coverage probability for the confidence intervals, we treat the random variable $\varsigma$ as a constant, since we do not have the original data from which it was estimated and since we do not know the data generating mechanism. This implies that in Figure 4 the term $S^{(\chi)}_{\hat{\mu}}$ in $(\ref{Sn})$ is not taken into account in the constriction of the confidence interval. 

Finally in Figure 5 we illustrate that even for a small number of repetitions (i.e.\ small effective sample sizes), we obtain good efficiency in the estimation of the wavelengths. 

\section{Conclusions}\label{secdisc}
The results here show that it is statistically possible to determine the neutron energy for a monochromatic beam with a good precision using multilayer neutron detectors. With relatively few layers ($ \leq 15$), already maximal information can be extracted and many layers do not significantly improve the precision of the results. 

For neutron beams with high intensity ($\lambda\geq 10^5$ particles), a statistical precision (width of 99 \% confidence interval) of less than $0.1$ \AA \, on the determination of the wavelength of the beam in the range 2.5-5 \AA \, is possible (Fig.\ref{muonn10}). Uncertainty in the neutron's cross section of the boron-10 isotope becomes dominant in the regime of high intensity beams and more than 10-20 layers. This means again that more than 10-20 layers are not needed (Fig.\ref{fci}).

An interesting further outcome of our work is that it shows that it might be possible, in high intensity experiments, with a precisely determined wavelength of a monochromatic neutron beam, to improve the statistical measurement of the boron-10 cross section by using an inverse of the method described in this manuscript. The systematic effects of such a measurement might be 
significant. In the limit of low intensity, a precision of 1 \AA \, in determining the wavelength of the monochromatic neutron beam is still possible. 

The asymptotic expansion used in the derivation of the asymptotic normality of the mle of the wavelength depends on two limit distribution results. The first is the asymptotic normality of the mle of the absorption probability $p$. Since we choose the effective number of neutrons that hits the detector ourselves, we are able to obtain an approximation which is as fine as wanted. Furthermore, the term $(\ref{Sp})$ in the total efficiency $(\ref{Sn})$, resulting from the mle of $p$, can be obtained as small as desired. A possible limitation here is that a large number of effective neutrons means running the experiment for a long time. In that case the assumption of a constant intensity Poisson process as a model may become questionable. A possible remedy for this is instead to do many repeated runs, while tightly controlling the experimental apparatus, in order to obtain a homogeneous Poisson process in each run. The second asymptotic result is the asymptotic normality of the estimator of $\varsigma$, which we conclude from Schmitt \textit{et al.} (1959). The number of data points used for the estimation of $\varsigma$ in that paper is 45, and therefore arguably on the boundary of what one can accept as an asymptotic normality result. A more serious practical limitation for us is that  we are not able to affect the term $(\ref{Shi})$ in  $(\ref{Sn})$ resulting from the estimator of  $\varsigma$. This puts a limit on the total efficiency that we can obtain for the wavelength estimation in our experimental setup. It also tells us, as noted above, that building a detector with many layers is not necessary, since for such a detector the term that we can affect in $(\ref{Sn})$ becomes negligible compared to term arising from the estimation of  $\varsigma$, and therefore increasing the number of layers will have negligible effect on $(\ref{Sn})$.

In a real detector there may be a degradation in the result achieved coming from systematic effects resulting from defects in the detector.

In this paper we have considered the Poisson process as a model for the incoming beam. Having real data it will in the future be possible to perform goodness of fit tests, e.g. for assessing the validity of the Poisson process model. A possible alternative model  for the incident beam is the negative binomial process. In fact, thinning of a negative binomial process also results in a negative binomial process, cf. Harremos \textit{et al.} (2007). However, unlike in the Poisson process case, the count processes $\{X_{i}(t)\}_{i\geq 1}$ in that case will not be independent, which makes the maximum likelihood approach more complicated. A possible solution could be to simplify the likelihood using some sort of quasi likelihood approach, e.g. by treating the count processes as independent and obtaining similar expressions for the likelihood as in this paper.  Model fit testing and negative binomial process modelling may be a direction for possible future research.


This manuscript concentrated on a monochromatic neutron beam. In the future our results will be generalised to discrete and continuous wavelength distributions for the incoming neutron beam.

\section{Acknowledgements}
VP's research is fully supported by the Swedish Research Council (SRC). The research of DA, RHW and KK is partially supported by the SRC. The authors gratefully acknowledge the SRC's support. The authors would furthermore like to thank the  associate editor and referees for their comments that have significantly improved the exposition and readability of the paper.
\\

Author's present address: {\"A}llingav{\"a}gen 12 lgh 1006 227 34 Lund \\
 
E-mail: pastuhov@maths.lth.se

\section{Appendix}
\noindent

\noindent
\textit{Proof Lemma \ref{Lmzero} \label{App2}.} 
For simplicity we skip the lower subscript $n$ but we assume that  $a$,  $b$, $c$, $d$ are as defined in (\ref{abcdx}). 

We study the monotonicity and convexity/concavity of $\tilde{f}$ on $[0, \infty )$ by studying the signs of $\tilde{f}'$ and $\tilde{f}''$ on $[0, \infty )$. For $k\geq2$ we have
\begin{eqnarray*}
\tilde{f}' &=& a(k+1)y^{k} - bky^{k-1} + c,\\
\tilde{f}'' &=& y^{k-2}k(a(k+1) y- b(k-1)).
\end{eqnarray*}

$(i):$ {\em The second derivative}.\\ 
Clearly $\tilde{f}''(0)=0$. Factoring out $k y^{k-2}\geq 0$, we see that to study the zeros and signs of $\tilde{f}''$ is equivalent to studying the zeros and signs of 
\begin{eqnarray*}
   g(y)&=&a(k+1)y - b(k-1),
\end{eqnarray*}
Clearly $g(0)=-b(k-1)<0$, $g(\infty) > 0$ and $g(y)$ has a unique root
\begin{eqnarray*}
y_{i.p.} = \frac{b(k-1)}{a(k+1)}.
\end{eqnarray*}
From the expressions in (\ref{abcdx}) we can see that both $a$ and $b$ are positive and $b > a$, which means that $y_{i.p.}\in(0, \infty)$.
  
Thus the function $\tilde{f}''$ is negative to the left of $y_{i.p.}$ and positive to the right of $y_{i.p.}$ which implies 
\begin{enumerate}
\item[a)] $\tilde{f}$ is concave on $(0,y_{i.p.})$, convex on $(y_{i.p.},\infty)$, and thus $y_{i.p.}$ is an inflection point for $\tilde{f}$.
\end{enumerate}

$(ii):$ {\em The first derivative}. We see that $\tilde{f}'(0)=c>0$. Furthermore
using the expressions for $a,b,c$ we see that $\tilde{f}'(1)=a(k+1)-kb+c   =0$. From the sign change of $\tilde{f}''$ at $y_{i.p.}$ we have that $\tilde{f}'$ is decreasing on $(0,y_{i.p.})$ and increasing on $(y_{i.p.},\infty)$. Now there are two possible cases:

${\bf Case\; A}:$ $y_{i.p.}<1$. In this case, the sign change of $\tilde{f}''$ together with $\tilde{f}'(0)=c>0, \tilde{f}'(1)=0$ and the continuity of $\tilde{f}$, implies that for some $y_1 < y_{i.p.}$, 
\begin{enumerate}
\item[b')] $\tilde{f}'$ is positive on $(0, y_1)$, negative on $(y_1, 1)$, positive on $(1,\infty)$,
\end{enumerate}
which of course implies 
\begin{enumerate}
\item[c')] $\tilde{f}$ is increasing on $(0, y_1)$, decreasing on $(y_1, 1)$, increasing on $(1,\infty)$.
\end{enumerate}

${\bf Case\; B}:$ $y_{i.p.}\geq 1$. In this case we know that $\tilde{f}'$ is decreasing and positive on $(0, 1)$, decreasing and negative on $(1, y_{i.p.})$ and increasing on $(y_{i.p.}, \infty)$. This implies that there is an $y_2$ such that $\tilde{f}'$ is negative on $(y_{i.p.}, y_2)$ and positive on $(y_2,\infty)$.
Thus the full statement becomes
\begin{enumerate}
\item[b'')] $\tilde{f}'$ is decreasing and positive on $(0,1)$, decreasing and negative on $(1, y_{i.p.})$, increasing and negative on $(y_{i.p.}, y_2)$, increasing and positive on $(y_2, \infty)$.
\end{enumerate}
which implies that
\begin{enumerate}
\item[c'')] $\tilde{f}$ is concave and increasing on $(0,1)$, concave and decreasing on $(1, y_{i.p.})$, convex and decreasing on $(y_{i.p.}, y_2)$, convex and increasing on $(y_2, \infty)$.
\end{enumerate}

$(iii):$ {\em The function}. We first note that $\tilde{f}(0)=-d<0$, and that the expression for the coefficients $a,b,c,d$ imply  $\tilde{f}(1)=a-b+c-d=0$. Now  we treat the two cases separately:

${\bf Case\; A}$: From the sign changes of $\tilde{f}''$ and $\tilde{f}'$, it follows that $\tilde{f}$ is concave and increasing on $(0,y_1)$, concave and decreasing on $(y_1, y_{i.p.})$, convex and decreasing on $(y_{i.p.}, 1)$. This together with $\tilde{f}(0)=-d<0$, $\tilde{f}(1)=0$ implies (and in fact only the information that $\tilde{f}$ is first increasing, then decreasing is enough) that there is a zero $\tilde{y}\in (0,1)$ for $\tilde{f}$.

${\bf Case\; B}$: In this case we have that $\tilde{f}$ is increasing and concave on $(0,1)$, which together with $\tilde{f}(0)=-d<0$, $\tilde{f}(1)=0$ implies that there is no zero for $\tilde{f}$ in the open $(0,1)$. 

Finally noting that a zero $\tilde{y}$ of $\tilde{f}$ in $(0,\infty)$, corresponds,
via $\tilde{y}=1-\tilde{p}$, to a zero $\tilde{p}$ of $f$ in $(-\infty,1)$, the Lemma follows. \eop

\noindent
\textit{Proof of Lemma \ref{Aoneroot} \label{App3}.}
From Lemma \ref{Lmzero}, we see that
\begin{eqnarray*}
 A_n = \left\{\frac{b_{n}(k - 1)}{a_{n}(k + 1)}<1\right\}.
\end{eqnarray*} 
We will prove that
\begin{eqnarray}\label{condrootml}
\frac{b_{n}(k - 1)}{a_{n}(k + 1)}\stackrel{a.s.}{\to} c,
\end{eqnarray} 
as $n\to\infty$, for some constant $c<1$. This immediately proves the condition of the lemma, since if $c<1$
\begin{eqnarray*}
\left\{\frac{b_{n}(k - 1)}{a_{n}(k + 1)}\to c \right\} \subseteq \underset{n\geq 1}{\cup}\underset{m\geq n}{\cap} A_{m}.
\end{eqnarray*} 

Now to prove (\ref{condrootml}), note that $\{s_j\}_{j=1}^n$ and $\{z_j\}_{j=1}^{n}$ in (\ref{abcdx}) are two sequences of i.i.d. random variables.  Thus from the strong law of large numbers 
\begin{eqnarray*}
	\frac{b_{n}(k - 1)}{a_{n}(k + 1)} &\stackrel{a.s}{\to}& \frac{k - 1}{k + 1}\frac{k - (k + 1)(1-p) +(1 - p)^{k-1}}{(k - 1) - k(1 - p) + (1 - p)^{k}} =: c,
\end{eqnarray*}
as $n\to\infty$. One can easily prove that $c<1$ by considering  the polynomial 
\begin{eqnarray*}
(k - 1)(1 - p)^{k + 1} - (k + 1)(1 - p)^{k} + (k + 1)(1 - p) - (k-1),
\end{eqnarray*}
which is negative for all $k > 1$ and $0 < p < 1$. This proves the lemma.
 \eop

\noindent
\textit{Proof of Theorem \ref{consnormef} \label{App4}.}

From Lemma \ref{Lmzero} it follows that there exists $n_{1}$ such that for all $n > n_{1}$ the mle $(\hat{p}_n, \hat{\lambda}_n)$ is a differentiable function of 
$(s_n, z_n)$, defined in (\ref{psilikk}). Therefore, the strong consistency of $(\hat{p}_n, \hat{\lambda}_n)$ follows from the strong law of large numbers and the continuous mapping theorem.

Next, $(s_n, z_n)$ is asymptotically normal, which follows from the central limit theorem. Using the delta method we prove the asymptotic normality of  $(\hat{p}_n, \hat{\lambda}_n)$.
\eop



\noindent
\textit{Proof of Corollary \ref{cormunorm}. }

Assume that there has been made $n$ measurements for $(\hat{p}_ {n},\hat{\lambda}_ {n})$ and $n'$ measurements for $\hat{\chi}_{n'}$, and that $(\hat{p}_ {n},\hat{\lambda}_ {n})$ and $\hat{\chi}_{n'}$ are independent. Let $n' = \lceil \gamma n \rceil $, with  $\gamma$ a proportionality factor that we introduce for convenience.


From the asymptotic normality of the estimators $\hat{p}_{n}$ and $\hat{\chi}_{n'}$ we have 
\begin{eqnarray}\label{normasymp}
\sqrt{n}(\hat{p}_{n} - p) \stackrel{d}{\to} \mathcal{N}(0,  \sigma^{2}_{p}),
\end{eqnarray} 
and
\begin{eqnarray}\label{normasymchi}
\sqrt{n}(\hat{\chi}_{n'} - \chi) = \sqrt{\frac{n}{n'}}  \sqrt{n'}(\hat{\chi}_{n'} - \chi) \nonumber \\
=\sqrt{\frac{n}{\lceil \gamma n \rceil}}  \sqrt{n'}(\hat{\chi}_{n'} - \chi) \stackrel{d}{\to} \mathcal{N}(0, \frac{\sigma^{2}_{\chi}}{\gamma}),
\end{eqnarray} 
as $n\to\infty$, since $\lim_{n \to \infty} \frac{n}{\lceil \gamma n \rceil} = \frac{1}{\gamma}$. Combining $(\ref{normasymp})$ and $\ref{normasymchi}$, the result follows from the delta method, see, for example, Chapter 3 in van der Vaart (1998). \eop

\end{document}